\documentclass[12pt, equation, leqno]{article}
\usepackage{amssymb}
\usepackage{amsmath}
\usepackage{theorem}
\evensidemargin\oddsidemargin
\makeatletter
\@addtoreset{equation}{section}

\makeatother
\newtheorem{theo}{Theorem}

\newtheorem{lem}{Lemma}
\newtheorem{cor}{Corollary}
{\theorembodyfont{\rmfamily}  \newtheorem{rem}{Remark}}
\newcommand{\PP}{\Bbb{P}}

\newcommand{\R}{\rm I\!R}

\newcommand{\E}{\Bbb{E}}

\newcommand{\hc}{\mathcal{H}}
\newcommand{\fc}{\mathcal{F}}
\newcommand{\be}{\begin{equation}}
\newcommand{\ee}{\end{equation}}
\newcommand{\nin}{\noindent}
\begin{document}

\date{}
\title{\ A generalization of Strassen's functional LIL}

\author{Uwe Einmahl\thanks{%
Research partially supported by an FWO-Vlaanderen Grant}\\ Departement Wiskunde\\
Vrije Universiteit Brussel\\ Pleinlaan 2 \\ B-1050 Brussel, Belgium\\
E-mail: ueinmahl@vub.ac.be}
\maketitle

\begin{abstract}
 \nin Let $ X_1, X_2, \ldots$ be a sequence of
 i.i.d. mean zero random variables and  let
$S_{n}$ denote the sum of the first $n$ random variables.   We show that whenever we have with probability one, $\limsup_{n \to \infty} |S_n|/c_n = \alpha_0 < \infty$ for a regular normalizing sequence $\{c_n\}$, the corresponding normalized partial sum process sequence  is relatively compact in $C[0,1]$ with canonical cluster set. Combining this result with some  LIL type results in the infinite variance case, we  obtain    Strassen type results in this setting. 
\end{abstract}

\bigskip\nin {\it Short title:} General Strassen-type results

\nin {\it AMS 2000 Subject Classifications:} 60F15, 60G50.

\nin {\it Keywords:} Hartman-Wintner LIL, Strassen- type results, infinite variance, very slowly varying functions,  sums of i.i.d. random variables, strong invariance principle. 

 \newpage

\section{Introduction}
Let $\{X,~X_{n};~n \geq 1 \}$ be a sequence 
of real-valued independent and identically distributed (i.i.d.)   
random variables, and let $S_{n} = \sum_{i = 1}^{n} X_{i}, ~~ n \geq 1$.
Define $Lx = \log_{e} \max \{e, x \}$ and  $LLx = L(Lx)$ for $x \in \R$.  
The classical Hartman-Wintner law of the iterated logarithm 
(LIL) states that
\begin{equation}
\limsup_{n \rightarrow \infty}
\pm S_{n}/(2nLL n)^{1/2}
= \sigma < \infty ~~\mbox{a.s.}
\end{equation} 
if and only if 
\begin{equation}
\E X = 0 ~~\mbox{and}~~\sigma^{2} = \E X^{2}  < \infty. \label{MO}
\end{equation}    

Starting with the work of Feller (1968) there has been quite some interest in finding extensions of the Hartman-Wintner LIL to the infinite variance case, especially for variables $X$ in the domain of attraction to the normal distribution. Recall that this means that one can find a centering sequence $\{\beta_n\} \subset \R$ and a normalizing sequence $\alpha_n \nearrow \infty$ so that
$$ \frac{S_n - \beta_n}{\alpha_n} \stackrel{d}{\to} Y \sim \mbox{ normal}(0,1),$$
where $\stackrel{d}{\to}$ stands for convergence in distribution.\\
Moreover, it is known that in this case, one has $\E |X|^p < \infty$ for all $0 < p <2$ and a possible choice for the centering sequence is given by $\beta_n = n\E X.$ The normalizing sequence $\alpha_n$ can be chosen of the form $\sqrt{ng(n)}$ where $g: [0,\infty[ \to [0,\infty[$ is  non-decreasing and slowly varying at infinity, that is we have $\lim_{t \to \infty} g(et)/g(t) =1.$\\
Under an additional symmetry assumption (which turned out to be unnecessary later) Feller (1968) proved in this case that 
\be
\limsup_{n \rightarrow \infty}
\pm S_{n}/\alpha_{[nLL n]}
= \sqrt{2}~~\mbox{a.s.}
\ee 
if and only if 
\begin{equation}
\E X = 0 ~~\mbox{and}~~\sum_{n=1}^{\infty}\PP\{|X| \ge \alpha_{[nLL n]}\}  < \infty. \label{FE}
\end{equation}
 It is easy to see that this result implies the Hartman-Wintner LIL since random variables satisfying condition (\ref{MO}) are trivially in the domain of attraction to the normal distribution with $\alpha_n = \sigma \sqrt{n}$ and $\beta_n = n\E X=0.$\\
Though this result is more general than the Hartman-Wintner LIL it was not the final word in the infinite variance case. The main problem is condition (\ref{FE}) which is not satisfied for many distributions of interest, most notably it cannot be applied for the distribution determined by the density $f(x)=1/|x|^3, |x| \ge 1.$ Another shortcoming is that it is essentially restricted to the domain of attraction case. This latter problem has been addressed by Klass (1976, 1977), Pruitt (1981), Li and Tomkins (2003) among other authors, but finding a normalizing sequence for the aforementioned example remained an open problem. \\ 
Einmahl and Li (2005) looked at  the following modification of the LIL problem:\\
Given a sequence, $a_n= \sqrt{nh(n)}$, where $h$ is a slowly varying non-decreasing function,  when does one have with probability one,
$$ 0 < \limsup_{n \to \infty} |S_n|/a_n < \infty \;?  $$
Somewhat surprisingly, it turned out that the Hartman-Wintner LIL could be generalized to a``law of the very slowly varying function'',  which is general enough to also give a normalizing sequence for the above example, that is for random variables with the distribution determined by the density $f(x)=1/|x|^3, |x| \ge 1.$\\
A natural question is now whether one can also find functional Strassen type versions of these results. All the results of this type known to us are restricted to the domain of attraction case (see, for instance, Kuelbs (1985), Einmahl (1989)).  We shall show  that whenever one has with probability one $\limsup_n |S_n|/c_n = \alpha_0 < \infty$ and the normalizing sequence $\{c_n\}$ is sufficiently regular, a functional LIL type result with canonical cluster set holds. Our result  implies among other things  that  the recent results of Einmahl and Li (2005) and also the two-sided version of the Klass (1976) LIL can be improved to functional LIL type results. 
\section{Statement of results}
To formulate our main result we need some extra notation.
 Let $S_{(n)}: \Omega \to C[0,1]$ the partial sum process of order $n$, that is,
$$
S_{(n)}(t)= S_{[nt]} + (nt -[nt])X_{[nt]+1}, 0 \le t\le 1.
$$
Strassen (1964) proved under condition (\ref{MO}) that with probability one, 
$$\{S_{(n)}/\sqrt{2nLLn}; n \ge 1\} \mbox{ is relatively compact in } C[0,1]$$ and
$$C(\{S_{(n)}/\sqrt{2nLLn}; n \ge 1\})= \sigma\mathcal{K},$$ 
where $C(\{x_n: n \ge 1\})$ denotes the  set of limit points of the sequence $\{x_n\} \subset C[0,1]$ which we shall call the ``cluster set'' of $\{x_n\}$ and
$$\mathcal{K}=\{f \in C[0,1]: f(t)=\int_0^t g(u)du, 0 \le t \le 1, \mbox{ where } \int_0^1 g^2(u)du \le 1\}.$$\\

\nin Recall that the Strassen LIL implies for any continuous functional $\psi: C[0,1] \to \R$:\\
\nin \emph{With probability one, $\{\psi(S_{(n)}/\sqrt{2nLLn}); n \ge 1\}$ is bounded and  the corresponding cluster set is equal to $\psi(\sigma \mathcal{K}).$}\\
\nin Choosing $\psi(f)=f(1), f \in C[0,1]$ one re-obtains the Hartman-Wintner LIL, but this is of course just one of many possible applications of the Strassen LIL.\\

\nin Let  $c_n$ be a sequence of real numbers satisfying the following two conditions,
\be
c_n/\sqrt{n} \nearrow \infty \label{RE}
\ee
and
\be
 \forall\; \epsilon >0\,\exists\, m_\epsilon \ge 1: c_n/c_m \le (1+\epsilon)(n/m), m_\epsilon \le m < n. \label{REG}
\ee  
Note that condition (\ref{REG}) is satisfied if 
 $c_n/ n$ is non-increasing or if
$c_n = c(n),$ where $c: [0,\infty) \to [0,\infty)$ is regularly varying at infinity with exponent $\gamma < 1$. \\
Let $H$ be the truncated second moment function of the random variable $X$, that is
$$H(t) = \E X^2 I\{|X| \le t\}, t \ge 0.$$
Set 
$$\alpha_0 = \sup\left\{\alpha \ge 0: \sum_{n=1}^{\infty} n^{-1}\exp\left(-\frac{\alpha^2 c^2_n}{2nH(c_n)}\right) = \infty\right\}
.$$
Then the main result of the present paper is as follows,
\begin{theo}\label{th1} Let $X, X_1, X_2, \ldots$ be i.i.d. mean zero random variables. 
Assume that 
\be
\sum_{n=1}^{\infty} \PP\{|X| \ge c_n\} < \infty, \label{MOM}
\ee
where $c_n$ is a  sequence of positive real numbers satisfying conditions (\ref{RE}) and (\ref{REG}).\\
If $\alpha_0 < \infty$ we have with probability one,
\be
\{S_{(n)}/c_n; n \ge 1 \} \mbox{ is relatively compact in }C[0,1] \label{COM}
\ee
and
\be
C(\{S_{(n)}/c_n; n \ge 1 \} ) = \alpha_0 \mathcal{K}. \label{CLU}
\ee
\end{theo}
\begin{rem} From relations (\ref{COM}) and (\ref{CLU}) it follows  that one has under assumption (\ref{MOM}) with probability one,
\be
\limsup_{n \to \infty} |S_n|/c_n = \alpha_0 < \infty. \label{eli}
\ee
Moreover we have $C(\{S_n/c_n; n \ge 1\})=[-\alpha_0, \alpha_0].$ Therefore, the above result implies Theorem 3 of Einmahl and Li (2005) if $\alpha_0 < \infty.$ It is easy to see that condition (\ref{MOM}) is necessary for (\ref{eli}) and thus also for (\ref{COM}) and (\ref{CLU}) to hold. 
\end{rem}
One of the main difficulties in applying (\ref{eli}) is the calculation of the parameter $\alpha_0$. A key result is Theorem 4, Einmahl and Li (2005) which makes it possible to calculate $\alpha_0$ (or finding at least bounds) in many cases of interest. Since we now know that the very same parameter determines the cluster set, we immediately have the functional versions of all LIL type results given in Einmahl and Li (2005) once  Theorem \ref{th1} has been proved.\\

\nin To formulate some of the possible corollaries, we need some extra notation. Let $\hc$ be the set of all continuous, non-decreasing functions $h: [0,\infty[ \,\to\, ]0,\infty[,$ which are slowly varying at infinity.\\
Set $f_{\tau} (t): = \exp ((Lt)^{\tau}), 0 \le \tau \le 1.$
Let $\hc_0 \subset  \hc$ be  the class of all functions so that 
$$
\lim_{t \to \infty} h(tf_{\tau}(t))/h(t) =1,
0 < \tau < 1.
$$
As this condition is more restrictive than slow variation, we call these functions  ``very slowly varying''.  Examples for functions in $\hc_0$ are
 $h(t)=(Lt)^r, r \ge 0$ and  $h(t)= (LLt)^p, p \ge 0.$ \\

\nin Combining Theorem 2, Einmahl and Li (2005) and Theorem \ref{th1} we obtain the following corollary which could be called the ``functional law of the very slowly varying function''. \\ Recall that $H$ is the truncated second moment function of the random variable $X.$
\begin{cor}
Let $X, X_1, X_2, \ldots$ be 
 i.i.d. random variables, and let
$S_{n} = \sum_{i=1}^{n} X_{i},~ n \geq 1$. Given a function $h \in \hc_0$ 
  set $\Psi(x) = \sqrt{xh(x)}$ and $a_{n} =\Psi(n), n \geq 1$. 
If there exists a constant $0 < \lambda < \infty$ such that     
\be
\E X = 0, ~~~\E \Psi^{-1}(|X|) < \infty 
\mbox{ and}~~~
\limsup_{x \rightarrow \infty} \frac{\Psi^{-1}(x LL x)}{x^{2} LLx}
H(x)   = \frac{\lambda^{2}}{2},   \label{COND}
\ee 
then we have with probability one,
 \be
 \{S_{(n)}/a_n; n \ge 1\} \mbox{ is relatively compact in }C[0,1] \label{COMP}
 \ee
 and
\be
C (\{S_{(n)}/a_n; n \ge 1\}) = \lambda \mathcal{K} \label{CLUSTER}
\ee  
Conversely, relations (\ref{COMP}) and (\ref{CLUSTER}) imply condition (\ref{COND}).  
\end{cor}
Similarly as in Einmahl and Li (2005) we can infer from the above result,
\begin{cor}
Let $p \geq 1$.  For any
constant $0 < \lambda < \infty$ we have with probability one,
$$
\{ S_{(n)}/ \sqrt{2n(LLn)^{p}}; n \ge 1\} \mbox{ is relatively compact in }C[0,1]
$$
and
$$
C (\{S_{(n)}/\sqrt{2n(LLn)^{p}}; n \ge 1\}) = \lambda \mathcal{K} 
$$ 
if and only if                                  
\begin{equation}
\E X = 0,\E X^{2}/(LL|X|)^{p} < \infty 
\mbox{ and}~
\limsup_{x \rightarrow \infty}(LLx)^{1-p}
H(x) = \lambda^{2}. \label{LOGLOG}
\end{equation}  
\end{cor}
\begin{rem}
If $p = 1$ then  condition (\ref{LOGLOG}) is equivalent to
\[
\E X = 0 ~~~\mbox{and}~~~\E X^{2} = \lambda^{2}.
\]
We  see that  the
Strassen LIL is a special case of Corollary 2.
\end{rem}
\begin{cor}
Let $r > 0$. For any
 constant $0 < \lambda < \infty$ we have with probability one,
$$
\{ S_{(n)}/ \sqrt{2n(Ln)^{r}}; n \ge 1\} \mbox{ is relatively compact in }C[0,1]
$$
and
$$
C (\{S_{(n)}/\sqrt{2n(Ln)^{r}}; n \ge 1\}) = \lambda \mathcal{K} 
$$ 
if and only if                                      
\begin{equation}
\E X = 0, ~\E X^{2}/(L|X|)^{r} < \infty 
\mbox{ and}~
\limsup_{x \rightarrow \infty} \frac{LLx}{(Lx)^{r}}
H(x) = 2^r \lambda^{2}.     \label{LOG}
\end{equation}
\end{cor}
We finally  show how the two-sided version of the Klass LIL can  be improved to a functional LIL. 
We need a certain function $K$ which is defined for any random variable $X: \Omega \to \R$ with $0 < \E|X| < \infty.$ Set 
$M(t): = \E |X|I\{|X|>t\}, t \ge 0.$
Then it is easy to see that the function 
$$G(t): = t^2/(H(t) + t M(t)), t >0$$
 is continuous and increasing and the function $K$ is defined
as its inverse function.  Moreover, one has for this function $K$ that as $x \nearrow \infty$
\be
K(x)/\sqrt{x} \nearrow \left(\E X^2\right)^{1/2} \in\, ]0,\infty]
\ee
and, 
\be
K(x)/x \searrow 0.
\ee
Set $\gamma_n = \sqrt{2}K(n/LLn)LLn, n \ge 1.$ 
Then $\gamma_n$ obviously satisfies conditions (\ref{RE}) and (\ref{REG}) and we can infer from our theorem,
\begin{cor}[Klass LIL]
Let $X,X_1,X_2,\ldots$ be iid mean zero random variables. Then we have for the partial sum process sequence $S_{(n)}$ based on these random variables,
\be
\{S_{(n)}/\gamma_n; n \ge 1\} \mbox{ is relatively compact in } C[0,1]
\ee
and
\be C(\{S_{(n)}/\gamma_n; n \ge 1\}) = \mathcal{K} \mbox{ a.s.}
\ee
if and only if
\be
\sum_{n=1}^{\infty} \PP\{|X| \ge \gamma_n\} < \infty \label{KL}
\ee
\end{cor}
Our proof of Theorem \ref{th1} is based on a new strong invariance principle for sums of i.i.d. random variables with possibly infinite variance which will be proved in Section 3. We shall show that up to an almost sure error term of order $o(c_n),$ the partial sum process sequence $S_{(n)}$ can be approximated by $\sigma_n W_{(n)}$ where $\sigma_n^2=H(c_n) \nearrow \E X^2 \in [0,\infty]$ and $W_{(n)}(t)=W(nt), 0 \le t \le 1$ with $\{W(s), s \ge 0\}$ being a standard Brownian motion. This way we can reduce the investigation of the asymptotic behavior of  $S_{(n)}/c_n$ to that one  of $\sigma_n W_{(n)}/c_n$ where some general results on Brownian motion due to Talagrand (1992) and Cs\'aki (1980), respectively, will come in handy.
\section{Proof of Theorem \ref{th1} }
Throughout the whole section we shall assume that $\{c_n\}$ is a sequence of positive real numbers
satisfying conditions (\ref{RE}) and  ($\ref{REG}$).
Moreover, $X, X_1,X_2,\ldots$ will always be a sequence of i.i.d. mean zero random variables satisfying
 \be
 \sum_{n=1}^{\infty}\PP\{|X| \ge c_n\} < \infty. \label{ass3}
 \ee
The following lemma collects a number of useful facts which we will need later on.
\begin{lem}
Assuming (\ref{ass3}) we have,
\begin{eqnarray} 
\sum_{n=1}^{\infty} \E |X|^3 I\{|X| \le c_n\}/c_n^3 < \infty \label{L1}\\
H(c_n)=\E X^2I\{|X| \le c_n\} = o(c^2_n/n) \mbox{ as } n \to \infty\label{L3} \\
M(c_n)=\E |X|I\{|X| > c_n\} = o(c_n/n) \mbox{ as }n \to \infty\label{L4}\\
\sum_{k=1}^n \E X I\{|X| \le c_k\} = o(c_n) \mbox{ as }n \to \infty \label{L6}
\end{eqnarray}
\end{lem}
{\bf Proof} The first three facts are already stated and proven in Lemma 1 of Einmahl and Li (2005).\\
The last fact follows by the same argument as on page 2037, Einmahl (1993). Just replace $\gamma_n$ by $c_n$ and use conditions (\ref{RE}) and (\ref{REG}).
$\Box$\\

\nin The crucial tool for our proof is the following infinite variance version of a strong invariance principle given in Einmahl and Mason (1993).\newpage
\begin{theo} Let $X$ be a mean zero random variable satisfying condition (\ref{ass3}). If the underlying p-space $(\Omega,\fc,\PP)$ is rich enough, one can define a sequence $\{X_n\}$ of independent copies of $X$ and a standard Brownian motion $\{W(t); t \ge 0\}$ such that we have as $n \to \infty$
\be \|S_{(n)} - \sigma_n W_{(n)}\|_{\infty} = o(c_n) \mbox{ a.s.}
\ee
where
$$\sigma_n^2 = H(c_n), n \ge 1$$
and $W_{(n)}(t) = W(nt), 0 \le t \le 1.$
\end{theo}
{\bf Proof} (STEP 1) We first note that using Theorem 3 of Einmahl (1988) in conjunction with fact (\ref{L6}), we can construct a sequence $\{X_n\}$ of independent copies of  $X$ and independent standard normal random variables $\{Y_n\}$ such that as $n \to \infty$
\be
\sum_{j=1}^n (X_j - \tau_j Y_j) = o(c_n) \mbox{ a.s.} \label{P1}
\ee
where $\tau_n^2 = \mathrm{Var}(XI\{|X| \le c_n\})\;, n \ge 1.$\\

\nin (STEP 2) We next claim that
\be
\sum_n \frac{\sigma_n - \tau_n}{c_n} Y_n  \mbox{  converges a.s.}\label{P2}
\ee
This implies via the Kronecker lemma  that  as $n \to \infty$
\be
\sum_{j=1}^n ({\sigma_j - \tau_j}) Y_j = o(c_n)\mbox{ a.s.} \label{P3}
\ee
To prove (\ref{P2}) we note that on account of $\E X=0,$
\begin{eqnarray*}
\sum_{n=1}^{\infty} \frac{(\sigma_n - \tau_n)^2}{c_n^2} &\le&\sum_{n=1}^{\infty} \frac{\sigma^2_n - \tau^2_n}{c_n^2}\\ &=&\sum_{n=1}^{\infty} \frac{(\E XI\{|X| > c_n\})^2}{c_n^2}\le \sum_{n=1}^{\infty} \frac{\epsilon_n^2}{n^2} < \infty
\end{eqnarray*}
where  $\epsilon_n = n\E |X|I\{|X| >c_n\}/c_n \to 0$ (see relation (\ref{L4})).\\

\nin (STEP 3) Set 
$$\Delta_n := \max_{1 \le k \le n} \left| \sum_{j=1}^k (\sigma_n - \sigma_j)Y_j \right|\;, n \ge 1
$$
We show that as $k \to \infty$
\be
\Delta_{2^k}/c_{2^k} \to 0 \mbox{ a.s.} \label{P4}
\ee
To that end we  note that by Kolmogorov's maximal inequality,
$$\PP\{\Delta_{2^\ell} \ge \epsilon c_{2^{\ell}}\} \le \epsilon^{-2}\sum_{k=1}^{2^{\ell} }\frac{(\sigma_{2^{\ell}} - \sigma_k)^2}{ c_{2^{\ell}}^2}\le \epsilon^{-2}\sum_{k=1}^{2^{\ell} }\frac{\sigma^2_{2^{\ell}} - \sigma_k^2}{c_{2^{\ell}}^2}
$$
Set $p_j = \PP\{c_{j-1} < |X| \le c_j\}, j \ge 1$ ($c_0=0$) Then we have:
\begin{eqnarray*}
\sum_{\ell =1}^{\infty} \sum_{k=1}^{2^{\ell} }\frac{\sigma^2_{2^{\ell}} - \sigma_k^2}{c_{2^{\ell}}^2}&\le& \sum_{\ell =1}^{\infty} \sum_{k=1}^{2^{\ell} }\sum_{j=k+1}^{2^{\ell}}p_j c_j^2/c_{2^{\ell}}^2\\
&=& \sum_{\ell =1}^{\infty}\sum_{j=2}^{2^{\ell}} (\sum_{k=1}^{j-1} 1)p_j c_j^2/c_{2^{\ell}}^2 \le \sum_{\ell =1}^{\infty}\sum_{j=2}^{2^{\ell}}jp_j c_j^2/c_{2^{\ell}}^2\\
&=&\sum_{j=2}^{\infty} jp_j c_j^2 \sum_{\ell: 2^{\ell} \ge j}c_{2^{\ell}}^{-2}
\le \sum_{j=2}^{\infty} j^2 p_j \sum_{\ell: 2^{\ell} \ge j}2^{-\ell} \le 2\sum_{j=2}^{\infty} jp_j < \infty
\end{eqnarray*}
Here we have used the fact that $c_j^2/c^2_{2^{\ell}}\le  j/2^{\ell}$ if $j \le 2^{\ell}$ which follows from condition (\ref{RE}).\\
We now obtain relation (\ref{P4}) via the Borel-Cantelli lemma.\\

\nin (STEP 4) Note that
$$\max_{2^{\ell} < n \le 2^{\ell +1}}\Delta_n \le \Delta_{2^{\ell +1}} + (\sigma_{2^{\ell +1}} - \sigma_{2^{\ell}})\max_{1 \le m \le 2^{\ell +1}}|T_m|,$$
where $T_m = \sum_{j=1}^m Y_j, m \ge 1.$\\
Moreover, we have by L\'evy's inequality,
\begin{eqnarray*}
\PP\left\{(\sigma_{2^{\ell +1}} - \sigma_{2^{\ell}})\max_{1 \le m \le 2^{\ell +1}}|T_m| \ge \epsilon c_{2^{\ell}}\right\}&\le& 2 \PP\left\{(\sigma_{2^{\ell +1}} - \sigma_{2^{\ell}})|T_{2^{\ell +1}}| \ge \epsilon c_{2^{\ell}}\right\}\\
&\le& 2\epsilon^{-2}(\sigma_{2^{\ell +1}} - \sigma_{2^{\ell}})^2 2^{\ell +1} c^{-2}_{2^{\ell}}\\ &\le&2\epsilon^{-2} (\sigma_{2^{\ell +1}}^2 - \sigma_{2^{\ell}}^2) 2^{\ell +1} c^{-2}_{2^{\ell }}
\end{eqnarray*}
Next observe that due to the monotonicity of $\sigma_n$ and in view of STEP 3  we have
$$
\sum_{\ell =1}^{\infty} 2^{\ell} \frac{\sigma_{2^{\ell +1}}^2 - \sigma_{2^{\ell}}^2}{c_{2^{\ell}}^2} \le \sum_{\ell =1}^{\infty} \sum_{k=1}^{2^{\ell } }\frac{\sigma^2_{2^{\ell +1}} - \sigma_k^2}{c_{2^{\ell +1}}^2} \cdot \frac{ c_{2^{\ell +1}}^2}{c_{2^{\ell }}^2}< \infty.
$$
Here we have used the fact that $c_{2^{\ell +1}}/c_{2^{\ell}} = O(1)$ due to assumption (\ref{REG}). We can infer via Borel-Cantelli that as $\ell \to \infty,$
\be
(\sigma_{2^{\ell +1}} - \sigma_{2^{\ell}})\max_{1 \le m \le 2^{\ell +1}}|T_m|
=o(c_{2^{\ell}}) \mbox{ a.s.} \label{P5}
\ee
In view of relation (\ref{P4}) this means that as $n \to \infty$
\be
\Delta_n/c_n \to 0 \mbox{ a.s.} \label{P6}
\ee
(STEP 5) Combining relations (\ref{P1}), (\ref{P3}) and (\ref{P6}) we conclude that as $n \to \infty$
$$\max_{1 \le k \le n} \left| \sum_{j=1}^k (X_j - \sigma_n Y_j) \right | = o(c_n) \mbox{ a.s.}
$$
Let $T_{(n)}$ be the partial sum process (of order $n$) determined by the sums $T_k, 1 \le k \le n.$ Then we can infer from the above relation that as $n \to \infty$
\be
\|S_{(n)} - \sigma_n T_{(n)}\|_{\infty} = o(c_n) \mbox{ a.s.} \label{P7}
\ee
W.l.o.g. we can assume that there exists a standard Brownian motion $\{W(t): t \ge 0\}$ satisfying $W(n) = T_n, n \ge 1$ which of course implies that
$$T_{(n)}(k/n) = W_{(n)}(k/n) = T_k, 1 \le k \le n.$$
By a standard argument (see, for instance, p. 485 of Einmahl and Mason (1993)),
it follows that
\be
\|T_{(n)} - W_{(n)}\|_{\infty} = O(\sqrt{\log n}) \mbox{ a.s.} \label{P8}
\ee
Recalling fact (\ref{L3}), we find that
\be
\|  \sigma_n T_{(n)} - \sigma_n W_{(n)}\|_{\infty} = O(\sigma_n \sqrt{\log n}) = o(c_n \sqrt{\log n /n}) = o(c_n) \mbox{ a.s.} \label{P9}
\ee
Combining relations (\ref{P7}) and (\ref{P9}) we obtain the assertion. $\Box$
\begin{rem}
After some  modification one can prove the above theorem more generally for sequences $\{c_n\}$ satisfying condition (\ref{REG}) and instead of condition (\ref{RE}):
$$\exists\, \alpha \ge 1/3: c_n/n^{\alpha} \mbox{ is noncreasing.}$$
Then this result also includes the strong invariance principle given in Einmahl and Mason (1993) so that the above theorem is in fact an extension of the earlier result. For our purposes the version above, however, is more than sufficient and it should be clear to the interested reader how to get the more general version.
\end{rem}
Due to invariance, it is now sufficient to prove the functional LIL type results for $\sigma_n W_{(n)}/c_n$. Here we proceed similarly as in the classical proof by Strassen (1964).\\
We first show for $\epsilon >0,$
\be
\PP\{\sigma_n W_{(n)}/c_n \in (\alpha_0\mathcal{K})^{\epsilon} \mbox{ eventually}\} =1. \label{outer}
\ee
This relation combined with the compactness of $\mathcal{K}$ already implies that with probability 1, $\{\sigma_n W_{(n)}/c_n; n \ge 1\}$ is relatively compact in $C[0,1].$ Moreover we can infer that with probability one,
\be
C(\{\sigma_n W_{(n)}/c_n; n \ge 1\}) \subset \alpha_0 \mathcal{K}.
\ee
Thus it remains to prove that with probability one,
\be
C(\{\sigma_n W_{(n)}/c_n; n \ge 1\}) \supset \alpha_0 \mathcal{K}.
\ee
Given that $C(\{\sigma_n W_{(n)}/c_n; n \ge 1\})$  is a closed subset of $C[0,1]$ as a cluster set, it suffices to show that we have for any function $f$ in the interior of $\alpha_0 \mathcal{K}$ and $\delta > 0,$ 
\be
\PP\{ \|f - \sigma_n W_{(n)}/c_n\|_{\infty} < \delta \mbox{ infinitely often}\}=1 \label{inner}
\ee
To prove (\ref{outer}) and (\ref{inner}) we need two further lemmas. We first show that the parameter $\alpha_0$ remains unchanged if we replace $H(c_n)$ by $H(\delta c_n),$ where $0 < \delta <1$. Though this also follows from Theorem 3, Einmahl and Li (2005), we prefer to give the direct (simple) proof since this makes our proof of Theorem \ref{th1} completely independent of the proofs given by Einmahl and Li (2005)
who used among other tools   a non-trivial result of Kesten (1970) on general cluster sets. 
\begin{lem}\label{lem1} Let $\alpha_0$ be defined as in Theorem \ref{th1}. Set
$$\tilde{\alpha}_0 = \sup\left\{\alpha \ge 0: \sum_{n=1}^{\infty} n^{-1}\exp\left(-\frac{\alpha^2 c^2_n}{2nH(\delta c_n)}\right) = \infty\right\},
$$
where $0 < \delta <1$.
Then we have: $\alpha_0 = \tilde{\alpha}_0.$
\end{lem}
{\bf Proof} Note that by monotonicity of $H$ we  obviously have: $\tilde{\alpha}_0 \le \alpha_0.$\\
Set $\Delta_n = H(c_n) - H(\delta c_n).$ 
Then we  trivially have $\Delta_n/c_n^2 \le \delta^{-1}\E |X|^3 I\{|X| \le c_n\}/c_n^3$ and fact (\ref{L1}) implies that
\be
\sum_{n=1}^{\infty} c_n^{-2} \Delta_n < \infty \label{e1}
\ee
Next observe that for any $0 < \epsilon <1,$
\begin{eqnarray*}
\frac{1}{n}\exp\left(-\frac{\alpha^2 c_n^2}{2nH(c_n)}\right) &=&
\frac{1}{n}\exp\left(-\frac{\alpha^2 c_n^2}{2n(H(\delta c_n) + \Delta_n)}\right)\\
&\le&\frac{1}{n}\exp\left(-\frac{\alpha^2 c_n^2}{2n(1+\epsilon)H(\delta c_n)}\right) +\frac{1}{n}\exp\left(-\frac{\alpha^2 c_n^2}{2n(1+\epsilon^{-1})\Delta_n}\right)\\
&\le& \frac{1}{n}\exp\left(-\frac{\alpha^2 c_n^2}{2n(1+\epsilon)H(\delta c_n)}\right) + \frac{2(1+\epsilon^{-1})}{\alpha^2} \cdot \frac{\Delta_n}{c_n^2}.
\end{eqnarray*}
In view of relation (\ref{e1}) this means that $\alpha_0 \le \sqrt{1+\epsilon}\;\tilde{\alpha}_0.$ Since we can choose $\epsilon$ arbitrarily small, we find that $\alpha_0 \le \tilde{\alpha}_0.$ $\Box$
\begin{lem}
Let $n_j \nearrow \infty$ a subsequence satisfying for large enough $j,$
$$ 1 < a_1 < n_{j+1}/n_j \le a_2 < \infty$$
Then we have:
\be
\sum_{j=1}^{\infty} \exp\left(-\frac{\alpha^2 c_{n_j}^2}{2n_j^2 \sigma_{n_j}^2}\right)  \begin{cases} =\infty & \mathrm{if}\;\alpha < \alpha_0\\ <\infty & \mathrm{if}\;\alpha > \alpha_0\end{cases} ,\label{geom}
\ee
where $\sigma_n^2 = H(c_n).$
\end{lem}
{\bf Proof.} 
Let  $\tilde{\sigma}_n^2 = H(\delta c_n),$ where $0 < \delta <1$ will be specified later.\\
Let $\alpha < \alpha_0 = \tilde{\alpha}_0.$ (See Lemma \ref{lem1}.) Since $c_n/\sqrt{n}$ is non-decreasing, we have
\begin{eqnarray*}
\infty &=& \sum_{j=j_0}^{\infty} \sum_{n=n_j +1}^{n_{j+1}}\frac{1}{n}\exp\left(-\frac{\alpha^2 c_n^2}{2n\tilde{\sigma}_n^2}\right)\\
&\le& (a_2 -1)\sum_{j=j_0}^{\infty} \exp\left(-\frac{\alpha^2 c_{n_j}^2}{2n_j\tilde{\sigma}_{n_{j+1}}^2}\right)\\
&\le&(a_2 -1)\sum_{j=j_0}^{\infty} \exp\left(-\frac{\alpha^2 c_{n_j}^2}{2n_j \sigma_{n_{j}}^2}\right)
\end{eqnarray*}
provided we have chosen $\delta$  small enough so that $\tilde{\sigma}_{n_{j+1}}^2=H(\delta c_{n_{j+1}}) \le H(c_{n_j})=\sigma^2_{n_j}.$ This is possible since $c_{n_{j+1}}/c_{n_j} = O(n_{j+1}/n_j) = O(1)$ (recall condition (\ref{REG})).
The other part of the lemma follows similarly. $\Box$\\

\nin {\bf Proof of (\ref{outer})} Our proof is based on the following inequality due to Talagrand (1992): Let $U$ be the unit ball in $(C[0,1],\|\cdot\|_{\infty})$. There exists an absolute constant $C > 0$ so that for all $t > 0$ and $\eta >0,$
\be
\PP\{W \not \in t\mathcal{K} + \eta U\} \le \exp\left(C\eta^{-2}-\frac{\eta t}{2} - \frac{t^2}{2}\right) \label{Tal}
\ee
Recall that $\sigma_n^2 = H(c_n).$ Using the fact that $W_{(n)} \stackrel{d}{=} \sqrt{n}W$, we get: 
\begin{eqnarray*}
\PP\{\sigma_n W_{(n)}/c_n \not \in \alpha_0\mathcal{K} + \epsilon U\}
&=& \PP\left\{W \not \in \frac{\alpha_0 c_n}{\sqrt{n\sigma_n^2}} \mathcal{K} + \frac{\epsilon c_n}{\sqrt{n\sigma_n^2}} U \right\}\\
&\le& \exp\left(C\frac{n\sigma_n^2}{\epsilon^2 c_n^2} - \frac{\alpha_0(\alpha_0 + \epsilon)c_n^2}{2n\sigma_n^2}\right)
\end{eqnarray*}
From fact  (\ref{L3}) we can infer that if $n$ is large enough, these probabilities are bounded above by
$$2  \exp\left(- \frac{\alpha_1^2 c_n^2}{2n\sigma_n^2}\right),$$
where $\alpha_1 = \sqrt{\alpha_0 (\alpha_0 + \epsilon)} > \alpha_0.$\newpage
\nin Let $\epsilon >0$ be fixed and set $n_j = [(1+\epsilon)^j], j \ge 1.$ Using  (\ref{geom}) in conjunction with the above bound we find that
$$
\sum_{j=1}^{\infty} \PP\{\sigma_{n_j} W_{(n_j)}/c_{n_j} \not \in \alpha_0\mathcal{K} + \epsilon U\} < \infty
$$
In view of the Borel-Cantelli lemma this means that 
$$\PP\{\sigma_{n_j} W_{(n_j)}/c_{n_j} \in \alpha_0\mathcal{K} + \epsilon U \mbox{ eventually}\}=1
$$
and consequently
\be
\limsup_{j \to \infty} d(\sigma_{n_j} W_{(n_j)}/c_{n_j}, \alpha_0\mathcal{K}) \le \epsilon \mbox{ a.s.}
\ee
where as usual $d(f,\alpha_0\mathcal{K})=\inf_{g \in \alpha_0\mathcal{K}}\|f-g\|_{\infty}, f \in C[0,1].$\\
Consider now $n_j \le n \le n_{j+1}$. From the definition of $\mathcal{K}$ it is clear that 
$$f \in \alpha_0\mathcal{K} \Rightarrow f(nt/n_{j+1}) \in \alpha_0 \mathcal{K}.$$ As we trivially have 
$$W_{(n)}(t)=W_{(n_{j+1})}(nt/n_{j+1}), 0 \le t \le 1$$
 we obtain for $n_j \le n \le n_{j+1}$ and large enough $j$
 \be 
 d\left(\frac{\sigma_n W_{(n)}}{c_n}, \frac{c_{n_{j+1}}}{c_n}\alpha_0 \mathcal{K}\right) \le 2\epsilon \mbox{ a.s.}
 \ee
 Recall that $\sigma_n$ is non-decreasing and use the fact that $t\mathcal{K} \subset \mathcal{K}, |t| \le 1.$\\
From assumption (\ref{REG}) and the definition of the subsequence $\{n_j\}$ we get for $n_j \le n \le n_{j+1}$ and large $j$, 
$$c_{n_{j+1}}/c_n \le c_{n_{j+1}}/c_{n_{j}} \le 1+2\epsilon.$$
Therefore, we have with probability one for large $j,$
\be
 d\left(\frac{\sigma_n W_{(n)}}{c_n}, (1+2\epsilon)\alpha_0\mathcal{K}\right) \le 2\epsilon, n_j \le n \le n_{j+1}.
\ee
Since $(1+a)\alpha_0\mathcal{K} \subset\alpha_0\mathcal{K} +a\alpha_0 U, a >0$, we can conclude that for large enough $n$
\be
 d\left(\frac{\sigma_n W_{(n)}}{c_n}, \alpha_0\mathcal{K}\right) \le 2(1+\alpha_0)\epsilon \mbox{ a.s.}
 \ee
and we have proven relation (\ref{outer}). $\Box$\\

\nin {\bf Proof of (\ref{inner})} To simplify notation, we set for $f \in C[0,1]$:
$$I(f)=\begin{cases} \left(\int_0^1 (f'(u))^2 du\right)^{1/2} &\mbox{if $f'$ exists}\\\infty &\mbox{otherwise}\end{cases}
$$
We use the following inequality given in Lemma 2 of Cs\'aki (1980).
\be
\PP\{\|W-f\|_{\infty} \le z\} \ge \exp(-I^2(f)/2)\PP\{\|W\|_{\infty} \le z\}, f\in C[0,1], z >0. \label{csa}
\ee
W.l.o.g. we shall assume that $0 < \alpha_0 <\infty.$\\
Let $\delta > 0$ be fixed. Set $n_j = m^j, j \ge 1,$ where $m \ge 2(1 +\alpha_0^2/\delta^2)$ is a natural number.\\
First note that if $j \ge 2$
\begin{eqnarray*}
\PP\left\{\sup_{0 \le t \le 1/m}\sigma_{n_j} |W(n_j t)| \ge \delta c_{n_j}\right\} &\le & 2\exp\left(-\frac{\delta^2 c_{n_j}^2}{2n_{j-1}\sigma^2_{n_j}}\right)\\&=& 2\exp\left(-\frac{m\delta^2 c_{n_j}^2}{2n_{j}\sigma^2_{n_j}}\right)
\le 2\exp\left(-\frac{2\alpha_0^2 c_{n_j}^2}{2n_{j}\sigma^2_{n_j}}\right),
\end{eqnarray*}
since  $m \ge 2\alpha_0^2/\delta^2.$
Recalling relation (\ref{geom}) we find that
\be
\sum_{j=2}^{\infty} \PP\left\{\sup_{0 \le t \le 1/m}\sigma_{n_j} |W(n_j t)| \ge \delta c_{n_j}\right\} < \infty \label{i1}
\ee
Therefore by the Borel-Cantelli lemma,
\be
\limsup_{j \to \infty} \frac{\sup_{0 \le t \le 1/m}\sigma_{n_j} |W(n_j t)| }{c_{n_j}} \le \delta \mbox{ a.s.} \label{i2}
\ee
Moreover, it is easy to see that
\be
\sup_{0 \le t \le 1/m} |f(t)| \le \alpha_0 /\sqrt{m} \le \delta, f \in \alpha_0 \mathcal{K}. \label{i3}
\ee
Combining relations (\ref{i2}) and (\ref{i3}), we can conclude for $f \in \alpha_0 \mathcal{K}$ that
\be
\liminf_{j \to \infty} \|f- \sigma_{n_j}W_{(n_j)}/c_{n_j}\|_{\infty}
\le 2\delta + \liminf_{j \to \infty} \sup_{1/m \le t \le 1}|f(t)- \sigma_{n_j}W(n_j t)/c_{n_j}|
\ee
which in turn (again by relation (\ref{i2})) is 
$$ \le 3\delta + \liminf_{j \to \infty} Z_j,$$
where
$$Z_j := \sup_{1/m \le t \le 1}|f(t)- \sigma_{n_j}\{W(n_j t)-W(n_{j-1})\}/c_{n_j}|, j \ge 2.$$
These variables are independent so that in view of the second Borel-Cantelli lemma we have $\liminf_{j \to \infty} Z_j \le 2\delta$ a.s.  and consequently
\be
\liminf_{j \to \infty} \|f- \sigma_{n_j}W_{(n_j)}/c_{n_j}\|_{\infty}
\le 5\delta \mbox{ a.s.} \label{i4}
\ee
if we can show,
\be
\sum_{j=2}^{\infty} \PP\{Z_j \le 2\delta\} = \infty. \label{i5}
\ee
It remains to prove (\ref{i5})  for functions $f$ in the interior of $\alpha_0 \mathcal{K}$ that is  under the assumption $I(f) < \alpha_0.$
 To that end we first note that
\begin{eqnarray*}
\PP\{Z_j \le 2\delta\} &\ge& \PP\left\{\sup_{1/m \le t \le 1}|f(t)- \sigma_{n_j}W(n_j t)/c_{n_j}| \le \delta, \sigma_{n_j}|W(n_{j-1})| \le \delta c_{n_j}\right\}\\
&\ge& \PP\{\|f-\sigma_{n_j}W_{(n_j)}/c_{n_j}\|_{\infty} < \delta\}
-\PP\{\sigma_{n_j}|W(n_{j-1})| \ge \delta c_{n_j}\}
\end{eqnarray*}
Since $\sum_j \PP\{\sigma_{n_j}|W(n_{j-1})| \ge \delta c_{n_j}\} < \infty$ (use (\ref{i1})), we see that (\ref{i5}) is proven once we have shown that
\be
I(f) < \alpha_0 \Rightarrow \sum_{j=1}^{\infty} \PP\{\|f-\sigma_{n_j}W_{(n_j)}/c_{n_j}\|_{\infty} < \delta\} =\infty \label{i6}
\ee
Applying inequality (\ref{csa}), we find that
\begin{eqnarray*}
\PP\{\|f-\sigma_{n_j}W_{(n_j)}/c_{n_j}\|_{\infty} < \delta\}&=&
\PP\{\|W - c_{n_j}f/(\sqrt{n_j}\sigma_{n_j})\|_{\infty} \le \delta c_{n_j}/(\sqrt{n_j}\sigma_{n_j})\}\\
&\ge& \exp\left(-\frac{I^2(f)c^2_{n_j}}{2n_j \sigma^2_{n_j}}\right) \PP\{\|W\|_{\infty} \le \delta c_{n_j}/(\sqrt{n_j}\sigma_{n_j})\}
\end{eqnarray*}
where we have used the trivial fact that $I(af) = a I(f), a >0.$\\
Recalling that $c_n/\sqrt{n\sigma_n^2} \to \infty$ (see (\ref{L3})), we can conclude that for large $j$
$$
\PP\{\|f-\sigma_{n_j}W_{(n_j)}/c_{n_j}\|_{\infty} < \delta\} \ge \frac{1}{2}\exp\left(-\frac{I^2(f)c^2_{n_j}}{2n_j \sigma^2_{n_j}}\right).
$$
In view of (\ref{geom}) it is now evident that (\ref{i6}) holds which in turn implies (\ref{i4}). We now see that relation (\ref{inner}) holds. $\Box$

\begin{center}
{\bf REFERENCES}
\end{center}
\begin{description}
\item
{\sc Cs\'aki, E.} (1980)
A relation between Chung's and Strassen's laws of the iterated logarithm.
{\em Z.Wahrsch. Verw. Gebiete} {\bf 54}~~ 287--301.
 \item 
 {\sc Einmahl, U.} (1988)  Strong approximations for partial sums of i.i.d. B-valued
r.v.'s in the domain of attraction of a Gaussian law. {\em  Probab. Theory
Related Fields} {\bf  77}~~ 65--85. 

\item
 {\sc Einmahl, U.} (1989) Stability results and strong invariance principles for
partial sums of Banach space valued random variables. {\it Ann.Probab}. {\bf 17}~~ 332--352.

\item
{\sc Einmahl, U.} (1993). Toward a general law of the iterated logarithm
in Banach space. {\em Ann. Probab.} {\bf 21}~~2012-2045.

\item
{\sc Einmahl, U.} and {\sc Li, D.} (2005) Some results on two-sided LIL behavior  {\it Ann. Probab.} {\bf 33}~~1601--1624.

\item
{\sc Einmahl, U.} and {\sc Mason, D.} (1993) Rates of clustering in Strassen's LIL for partial sum processes 
{\it Probab. Th. Rel. Fields} {\bf 97}~~479--487.

\item
{\sc Feller, W.} (1968). An extension of the law of the iterated logarithm
to variables without variance. {\em J. Math. Mech.} {\bf 18}~~343-355.

\item
{\sc Kesten, H.} (1970). The limit points of a normalized random walk.
{\em Ann. Math. Statist.} {\bf 41}~~1173-1205.

\item
{\sc Klass, M.} (1976). Toward a universal law of the iterated logarithm  
I.  {\em Z. Wahrsch. Verw. Gebiete} {\bf 36}~~165-178.

\item
{\sc Klass, M.} (1977). Toward a universal law of the iterated logarithm
II.  {\em Z. Wahrsch. Verw. Gebiete} {\bf 39}~~151-165.

\item
{\sc Kuelbs, J.} (1985). The LIL when $X$ is in the domain of attraction of a Gaussian law. {\em Ann. Probab.} {\bf 13}~~825--859.

\item
{\sc Li, D. L.} and {\sc Tomkins, R. J.} (2003). The law of the logarithm for weighted sums of independent random variables. {\em J. Theoretical Probab.}~{\bf 16}~~519--542.
 
\item
{\sc Pruitt, W.} (1981). General one-sided laws of the iterated logarithm.
{\em Ann. Probab.} {\bf 9}~~1-48.

\item
{\sc Strassen, V.} (1964). An invariance principle for 
the law of the iterated logarithm. {\em Z. Wahrsch. Verw. Geb.}~{\bf 3}
~~211-226.

\item
{\sc Talagrand, M.} (1992). On the rate of clustering in Strassen's LIL for Brownian motion. In: \emph{Probability in Banach spaces 8} (R. Dudley, M. Hahn, J. Kuelbs, eds.) 339--347. Birkh\"auser, Boston.

\end{description}
\end{document}